\documentclass[11pt,a4paper]{amsart}

\usepackage[utf8x]{inputenc} 
\usepackage[a4paper]{geometry} 
\usepackage{enumerate} 
\usepackage{amsmath} 
\usepackage{amssymb}
\usepackage{color,colortbl}
\usepackage[dvipsnames,table]{xcolor}
\usepackage[most]{tcolorbox}
\usepackage{fancyvrb}   
 
\definecolor{airforceblue}{rgb}{0.36, 0.54, 0.66}
\definecolor{bleudefrance}{rgb}{0.19, 0.55, 0.91}
\definecolor{darkorchid}{rgb}{0.6, 0.2, 0.8}
\definecolor{darkorange}{rgb}{1.0, 0.55, 0.0}
\definecolor{darkspringgreen}{rgb}{0.09, 0.45, 0.27}
\definecolor{commentoutput}{rgb}{0.40, 0.00, 0.0}
\definecolor{output}{rgb}{0.8, 0.0, 0.0}  
\definecolor{circOut}{rgb}{0.4, 1.0, 0.0} 

\definecolor{Gray}{gray}{0.7}  

\usepackage[arrow,curve,matrix,tips,frame,color]{xy}
\usepackage{adjustbox}
\usepackage{multirow}

\usepackage{verbatim}
\usepackage{hyperref}
                 \hypersetup{ pdfborder={0 0 0}, 
                              colorlinks=true, 
                              linktoc=page,
                              pdfauthor={Giovanni Staglian{\`o}}, 
                              pdftitle={Computations with rational maps between multi-projective varieties}
                            }                         

\theoremstyle{plain} 
\newtheorem{proposition}{Proposition}[section] 
\newtheorem{prop-def}[proposition]{Proposition-Definition}

\theoremstyle{definition}

\newtheorem{example}[proposition]{Example} 
\theoremstyle{remark} 
\newtheorem{remark}[proposition]{Remark}

\newcommand{\ZZ}{{\mathbb{Z}}} 
         
\newcommand{\PP}{{\mathbb{P}}}

\providecommand{\sat}{\mathop{\rm sat}\nolimits} 
\providecommand{\Dom}{\mathop{Dom}\nolimits} 
\providecommand{\codim}{\mathop{\rm codim}\nolimits}

\numberwithin{equation}{section}

\title{Computations with rational maps between multi-projective varieties} 
\address{Dipartimento di Matematica e Informatica, Universit\`a degli Studi di Catania} 
\author[G. Staglian\`o]{Giovanni Staglian\`o}
\email{\href{mailto:giovannistagliano@gmail.com}{giovannistagliano@gmail.com}} 
\subjclass[2010]{
     14E05, 
     14Q15
}

\begin{document}

\begin{abstract} 
We briefly describe the algorithms behind some of the functions provided by the \emph{Macaulay2} package 
\href{https://faculty.math.illinois.edu/Macaulay2/doc/Macaulay2/share/doc/Macaulay2/MultiprojectiveVarieties/html/index.html}{\emph{MultiprojectiveVarieties}},
a package for multi-projective varieties and 
rational maps between them.
\end{abstract}

\maketitle

\section*{Introduction}
This paper is a natural sequel of \cite{packageCremona},
where we presented some of the algorithms 
implemented in the \emph{Macaulay2} package \href{https://faculty.math.illinois.edu/Macaulay2/doc/Macaulay2/share/doc/Macaulay2/Cremona/html/index.html}{\emph{Cremona}} \cite{CremonaPackageSource},
related to computations with rational and birational maps 
between closed subvarieties of projective spaces.

Here 
we describe methods for working with 
rational and birational maps between multi-projective varieties, that is,
closed subvarieties of products of projective spaces.
For instance, we explain how to compute the degrees of such maps, 
their graphs, 
and the inverses when they exist.
All these methods are implemented in the
\emph{Macaulay2} package \href{https://faculty.math.illinois.edu/Macaulay2/doc/Macaulay2/share/doc/Macaulay2/MultiprojectiveVarieties/html/index.html}{\emph{MultiprojectiveVarieties}}.%

From a theoretical point of view, we know that every multi-projective variety 
is isomorphic, via the Segre embedding, to 
a projective variety embedded into a single projective space. Therefore,
every rational map 
between multi-projective varieties
can be regarded as a 
 rational map
between ordinary subvarieties of projective spaces.
This, however, introduces a lot of new variables, making computation more difficult.

Moreover,
 basic constructions on rational maps
naturally lead one to consider rational maps 
between multi-projective varieties.
For instance,  the graph 
of a rational map is a
closed subvariety of the 
product of the source and of the target of the map. 
Using 
the package \emph{Cremona},
it is generally easy to verify that the first projection from the graph 
is birational, 
but to calculate, for instance, its inverse 
we need the tools provided by the package presented here.

In section~\ref{section 1}, 
we give a concise overview of the theory of rational maps between multi-projective varieties,
emphasizing the computational aspects
and making clear how they can be
 represented in a computer.
 For more details on the theory see, \emph{e.g.}, \cite{harris-firstcourse,hartshorne-ag}. 
In section~\ref{section 2},
with the help of an example, we show 
how one can work with such maps using \emph{Macaulay2} \cite{macaulay2}.

 \section{An overview of rational maps between multi-projective varieties}\label{section 1}

\subsection{Notation and terminology}
 Throughout this paper, we keep the following notation.
 Let $K$ denote an arbitrary field. 
 Consider
  the following polynomial ring in $r$ groups of variables
 \begin{equation*} R = K[x_0^{(1)},\ldots,x_{n_1}^{(1)};\ldots ; x_0^{(r)},\ldots,x_{n_r}^{(r)}],\end{equation*}
 equipped with the  
 $\mathbb{Z}^r$-grading, where the degree of each variable is a standard basis vector.
 More 
 precisely, we set 
 $\deg(x_i^{(j)}) = (0,\ldots,0,1,0,\ldots,0)\in\mathbb{Z}^r$, where $1$ occurs at position $j$; we call this the \emph{standard 
 $\mathbb{Z}^r$-grading} on $R$.
 The polynomial ring $R$ is the homogeneous coordinate ring of 
 the product of $r$ projective spaces 
 \begin{equation*} \mathbf{P}^{n_1,\ldots,n_r}= \PP^{n_1}\times \cdots \times \PP^{n_r}.\end{equation*}
 The closed subsets (of 
 the Zariski topology) of
 $\mathbf P^{n_1,\ldots,n_r}$ are of the form 
 \[V(\mathfrak a) = \{p\in\mathbf P^{n_1,\ldots,n_r}: F(p)=0\mbox{ for all homogeneous }F\in \mathfrak a\},\]
where $\mathfrak a$ is a homogeneous ideal in $R$.
 For any homogeneous ideal $\mathfrak a\subseteq R$, the \emph{multi-saturation} of $\mathfrak a$
 is the homogeneous ideal
 \[
 \sat(\mathfrak a) = \left(\cdots\left(\left(\mathfrak a:(x_0^{(1)},\ldots,x_{n_1}^{(1)})^{\infty}\right):(x_0^{(2)},\ldots,x_{n_2}^{(2)})^{\infty}\right):\cdots\right):(x_0^{(r)},\ldots,x_{n_r}^{(r)})^{\infty} .
 \]
 One says that $\mathfrak a$ is \emph{multi-saturated}  
 if $\mathfrak a = \sat (\mathfrak a)$.
 Two homogeneous ideals $\mathfrak a,\mathfrak a'\subseteq R$ define the same subscheme of $\mathbf P^{n_1,\ldots,n_r}$
   if and only if ${\sat(\mathfrak a)} = {\sat(\mathfrak a'})$,
   and they define the same subset if and only if
   $\sqrt{\sat(\mathfrak a)} = \sqrt{\sat(\mathfrak a')}$.
   
 We fix a homogeneous absolutely prime ideal $I\subset R$,
 and we may also assume that
 $I$ is multi-saturated.
 The graded domain 
 $R/I $
 is the homogeneous coordinate ring  
 of an absolutely irreducible multi-projective variety 
 \[X=V(I)\subseteq \mathbf{P}^{n_1,\ldots,n_r}= \PP^{n_1}\times \cdots \times \PP^{n_r} .\]
 There is a similar correspondence between homogeneous ideals in $R/I$ 
 and closed subsets of $X$. 
 The two most important invariants of $X$ are:
the dimension  (as a topological space), which 
   is 
   the (Krull) dimension
  of the homogeneous coordinate ring $R/I$ minus $r$; and 
 the multi-degree, an integral homogeneous polynomial 
 of degree $\codim X = n_1+\cdots+n_r-\dim X$
 in $r$ variables (see \cite[Lecture~19]{harris-firstcourse} and \cite[p.~165]{millersturmfels}).
   
 Similarly, let us take another polynomial ring 
 in $s$ groups of variables
 \begin{equation*} S = K[y_0^{(1)},\ldots,y_{m_1}^{(1)};\ldots ; y_0^{(s)},\ldots,y_{m_s}^{(s)}],\end{equation*}
 equipped with the standard  
 $\mathbb{Z}^s$-grading.
 Let $J\subset S$ be a multi-saturated homogeneous absolutely prime ideal,
 and let \[Y=V(J)\subseteq \mathbf P^{m_1,\ldots,m_s} = \PP^{m_1}\times\cdots\times\PP^{m_s}\] be 
 the absolutely irreducible multi-projective variety defined by $J$.

 \subsection{Rational maps to an embedded projective variety}
 In this subsection we consider the particular case when $s=1$,
 and we set $\PP^m=\mathbf P^{m_1,\ldots,m_s}$. Then $Y\subseteq \PP^m$
 is an embedded projective variety.
  \subsubsection{Definition of rational map}
We call \emph{multi-form} (or simply \emph{form})  a homogeneous element of $R/I$.
 To a vector $\mathbf{F}=(F_0,\ldots,F_m)$ of $m+1$ forms in $R/I$ 
 of the same multi-degree, which are not all zero, 
 we associate a continuous map \[\phi_{\mathbf F}:X\setminus V(\mathbf F)\longrightarrow \PP^m, 
 \mbox{ defined by }
 p\in X\setminus V(\mathbf F)  \stackrel{\phi_{\mathbf F}}{\longmapsto} (F_0(p),\ldots,F_m(p))\in\PP^m .\]
 If $\mathbf G = (G_0,\ldots,G_m)$ 
 is another such vector of forms in $R/I$ of the same multi-degree,
 then we say that 
 $\mathbf F \sim \mathbf G$ 
 if 
 $\phi_{\mathbf F}(p) = \phi_{\mathbf G}(p)$ for each  
 $p\in X\setminus (V(\mathbf F)\cup V(\mathbf G))$.
   We 
   have $\mathbf F \sim \mathbf G$ if and only if 
  $\phi_{\mathbf F} = \phi_{\mathbf G}$
 on some nonempty open subset $U$ of $X\setminus (V(\mathbf F)\cup V(\mathbf G))$;
 in particular $\sim$ is an equivalence relation. 
  A \emph{rational map} $\Phi:X\dashrightarrow Y$ 
  is defined as an equivalence class 
  of non-zero vectors of $m+1$ forms $\mathbf F=(F_0,\ldots,F_m)$ in $R/I$
  of the same multi-degree,
  with respect to the relation $\sim$, such that  
  for same (and hence  every) representative $\mathbf F$ 
  we have that the image of 
  $\phi_{\mathbf F}$ is contained in $Y$.
  If $p\in X\setminus V(\mathbf F)$ for some representative $\mathbf F$, 
  we set $\Phi(p) = \phi_{\mathbf F}(p)$ and we say that $\Phi$ is \emph{defined} at $p$.
  The \emph{domain} of $\Phi$, denoted by $\Dom(\Phi)$,
  is the set of points where $\Phi$ is defined, that is, it
  is the largest open subset of $ X$ such that 
  the map $\phi_{\mathbf{  F}}$ is defined for 
  some representative 
  $\mathbf F$.  
  The complementary set in $X$ of the domain of $\Phi$
  is called \emph{base locus}.
   A rational map $\Phi:X\dashrightarrow Y$ is called a \emph{morphism}
  if it everywhere defined, that is, 
if its base locus is empty.

\subsubsection{Establishing the equality of rational maps}
 Notice that if a vector $\mathbf F = (F_0,\ldots,F_m)$ of forms
 in $R/I$ 
 represents a rational map $\Phi:X\dashrightarrow Y$, then also the vector  
 $H\cdot \mathbf F = (H F_0,\ldots,H F_m)$ represents $\Phi$,
 for each nonzero form  $H$ in $R/I$.
 More generally, 
  two vectors  $\mathbf F = (F_0,\ldots,F_m)$ and 
  $\mathbf G = (G_0,\ldots,G_m)$, as the ones considered above,
  represent the same rational map $\Phi:X\dashrightarrow Y$
  if and only if 
  \[\mathrm{rk} \begin{pmatrix} F_0 & \cdots & F_m \\ G_0 & \cdots & G_m \end{pmatrix} < 2 ,\]
  that is, if and only if $F_i G_j - F_j G_i$ vanishes identically on $X$,
  for every $i,j=0,\ldots,m$.

\subsubsection{Determining the domain of a rational map}
 Let $\Phi:X \dashrightarrow Y$
 be a rational map and let $\mathbf F = (F_0,\ldots, F_m)$ be
 one of its representatives.
 A \emph{syzygy} of $\mathbf F$ 
 is a vector $\mathbf H = (H_0,\ldots, H_m)$ of forms  in $R/I$
 such that $\sum_{i=0}^m H_i F_i = 0$.
 Let $M_{\mathbf F}$ be a matrix 
 whose columns form a  set of generators 
 for the module of syzygies of $\mathbf F$.
 The following result is proved in 
 \cite[Proposition 1.1]{Simis2004162}, although stated there only for $r=1$.
 \begin{proposition}
  The representatives of the rational map $\Phi$
  correspond bijectively to the 
  homogeneous vectors in the rank one graded $(R/I)$-module 
  \[ \ker(M_{\mathbf F}^{t}) \subset (R/I)^{m+1} .\]
  \end{proposition}
Let $\mathbf F_1,\ldots,\mathbf F_p$ be a set 
  of minimal homogeneous generators of $\ker(M_{\mathbf F}^{t})$. 
  The base locus of $\Phi$
  is the closed subset of $X$ 
   where all the entries of $\mathbf F_i$, for $i=1,\ldots,p$, vanish.
    The sequence of multi-degrees $(\deg \mathbf F_1 , \ldots , \deg \mathbf F_p)$, defined up to ordering,
  is called 
  the \emph{degree sequence} of $\Phi$.
  \begin{example}
  In the case when
  $R/I$  is  a unique factorization domain (\emph{e.g.}, $X=\PP^{n_1}\times\cdots\times\PP^{n_r}$), then
  a rational map $\Phi:X\dashrightarrow Y$ is uniquely represented up to proportionality, that is,
  the degree sequence of $\Phi$ consists of a unique element. 
\end{example}

 \subsubsection{Direct and inverse images via rational maps}
 Let $\Phi:X\dashrightarrow Y$ be a rational map,
 and let $\mathcal M$ be a set of generators for the $(R/I)$-module of
 representatives of $\Phi$.  For  
 $\mathbf F = (F_0,\ldots,F_m)\in \mathcal M$,
 we consider the graded $K$-algebra homomorphism $\varphi_{\mathbf F}:S/J\to R/I$ 
 defined by $\varphi_{\mathbf F}(y_i) = F_i\in R/I$.
 
 For each homogeneous ideal $\mathfrak a \subseteq R/I$ (resp. $\mathfrak b\subseteq S/J$),
 we have a closed subset $V(\mathfrak a)\subseteq X$ (resp. $V(\mathfrak b)\subseteq Y)$.
 The \emph{direct image}  of $V(\mathfrak a)$ via $\Phi$, denoted by $\overline{\Phi(V(\mathfrak a))}$,
 and 
 the \emph{inverse image}  of $V(\mathfrak b)$ via $\Phi$, denoted by $\overline{\Phi^{-1}(V(\mathfrak b))}$,
 as sets, are given by the following closure:
 \[
 \overline{\Phi(V(\mathfrak a))} =\overline{\{\Phi(p):p\in\Dom(\Phi)\cap V(\mathfrak a)\}},\quad 
 \overline{\Phi^{-1}(V(\mathfrak b))} =\overline{\{p\in\Dom(\Phi):\Phi(p)\in V(\mathfrak b)\}}.
 \]
 The following 
 result follows from elementary commutative algebra, and
 it tells us how to calculate
 direct and inverse images.
\begin{proposition}\label{direct and inverse image} The following formulas hold:
  \begin{align*}
   \overline{\Phi(V(\mathfrak a))} = & 
    \bigcup_{\mathbf F\in \mathcal M} V\left(\varphi_{\mathbf F}^{-1}(\mathfrak a)\right) 
   = V\left(\bigcap_{\mathbf F\in \mathcal M}\varphi_{\mathbf F}^{-1}(\mathfrak a)\right) ;
   \\
   \overline{\Phi^{-1}(V(\mathfrak b))} 
=& 
 \bigcup_{\mathbf F\in \mathcal M} V\left(\varphi_{\mathbf F}(\mathfrak b):(\mathbf F)^{\infty}\right)
=V\left(\bigcap_{\mathbf F\in \mathcal M} \varphi_{\mathbf F}(\mathfrak b):(\mathbf F)^{\infty}\right) .
\end{align*}
\end{proposition}
As a consequence, we obtain that
if $\mathbf F$ is any of the representatives of $\Phi$,
then
\[\overline{\Phi(X)}=V(\ker \varphi_{\mathbf F}) .\]
The direct image $\overline{\Phi(X)}$
is called the (closure of the) \emph{image} of $\Phi$.
We say that $\Phi$ is  \emph{dominant} if $\overline{\Phi(X)}=Y$.

 \subsection{Rational maps to a multi-projective variety}
We now consider the general case when $s\geq 1$,
and hence $Y\subseteq \mathbf P^{m_1,\ldots,m_s}= \PP^{m_1}\times\cdots\times\PP^{m_s}$
is a multi-projective variety.
Let us denote by $\pi_i:\mathbf P^{m_1,\ldots,m_s}\to \PP^{m_i}$ 
the $i$-th projection, and let $Y_i=\pi_i(Y)$.
\subsubsection{Definition of multi-rational map}
We define a \emph{multi-rational map} (or simply rational map) 
\[\Phi:X\dashrightarrow Y\]
as an $s$-tuple of rational maps $\Phi_i:X\dashrightarrow \PP^{m_i}$
such that the image of $\Phi_i$ 
is contained in $Y_i$, for $i=1,\ldots,s$.
The domain of a multi-rational map $\Phi$ is the 
intersection 
\[
 \Dom(\Phi) = \bigcap_{i=1}^s \Dom(\Phi_i) .
\]
In other words, $\Phi$ is defined at a point $p\in X$ if and only if $\Phi_i$ is defined at $p$ for all $i=1,\ldots,s$,
and in that case  we set $\Phi(p)=(\Phi_1(p),\ldots,\Phi_s(p))\in\mathbf{P}^{m_1,\ldots,m_s}$.
Analogously with the case $s=1$,
we call \emph{base locus} the complementary set in $X$
of the domain of $\Phi$, and we say that $\Phi$ is a 
 \emph{morphism}
  if 
  $X = \Dom(\Phi)$.
We say that $\Phi$ is \emph{dominant} if for some (and hence  every)
  open subset $U$ of the domain of $\Phi$, the set $\{\Phi(p):p\in U\}$ is dense in $Y$.
     
\subsubsection{Composition of multi-rational maps}
If $\Psi=(\Psi_1,\ldots,\Psi_t):Y\dashrightarrow Z$ 
is another multi-rational map, then $\Phi$ and $\Psi$
 can be composed if $\Phi(\Dom(\Phi))\cap \Dom(\Psi)\neq \emptyset$;
 in particular, this happens when either $\Phi$ is dominant or $\Psi$ is a morphism.
  If $\mathbf F^{(1)},\ldots,\mathbf F^{(s)}$ are, respectively, representatives of 
 $\Phi_1,\ldots,\Phi_s$, and if $\mathbf G^{(j)}$ is a representative of $\Psi_j$,
 then the vector $\mathbf G^{(j)}(\mathbf F^{(1)},\ldots,\mathbf F^{(s)})$ 
 is a representative of $(\Psi\circ \Phi)_j=\Psi_j\circ\Phi$.
 
  So we can consider the category 
  of (multi)-projective varieties 
  and dominant (multi)-rational maps. An ``\emph{isomorphism}'' in this category 
  is called a birational map, that is, 
   $\Phi:X\dashrightarrow Y$ is a birational map 
   if it admits an inverse, namely a multi-rational map $\Phi^{-1}:Y\dashrightarrow X$
   such that $\Phi^{-1} \circ \Phi = \mathrm{id}_X$ 
   and $\Phi\circ \Phi^{-1} = \mathrm{id}_Y$ as (multi)-rational maps.
  A birational morphism $\Phi:X\dashrightarrow Y$ is called \emph{isomorphism} 
  if  $\Phi^{-1}$ is a morphism.
Also (multi)-projective varieties and morphisms form a category.

\subsubsection{Example: the Segre embedding}
 Let 
 $N=(n_1+1)\cdots(n_r+1)-1$,
 and let us
 consider $\PP^N$ 
 with the homogeneous coordinate ring 
 $
 K[z_{(\iota_1,\ldots,\iota_r)}:\, \iota_j=0,\ldots,n_j,\, j=1,\ldots,r]
 $,
 where the variables are the entries
 of the generic $r$-dimensional matrix of shape $(n_1+1)\times\cdots\times (n_r+1)$.
 The \emph{Segre embedding} of $\PP^{n_1}\times \cdots \times \PP^{n_r}$ into $\PP^N$
 is the rational map
 \begin{gather*}
 \mathfrak S_{n_1,\ldots,n_r}:\PP^{n_1}\times \cdots \times \PP^{n_r}\dashrightarrow \PP^N, \\
 \end{gather*}
represented  by the following ring map:
 \[
 \begin{gathered} 
 K[z_{(\iota_1,\ldots,\iota_r)}:\, \iota_j=0,\ldots,n_j,\, j=1,\ldots,r] \to K[x_0^{(1)},\ldots,x_{n_1}^{(1)},\ldots,x_0^{(r)},\ldots,x_{n_r}^{(r)}], \\
 z_{(\iota_1,\ldots,\iota_r)} \mapsto x_{\iota_1}^{(1)}\cdots x_{\iota_r}^{(r)} .
 \end{gathered}
\]
This ring map (or better the forms defining it) represents uniquely up to proportionality
 the rational map $\mathfrak S_{n_1,\ldots,n_r}$, and it is also clear that 
it is an injective morphism.
The image of $\mathfrak S_{n_1,\ldots,n_r}$ is 
the projective variety of all $r$-dimensional matrices of rank $1$.
If we consider $\mathfrak S_{n_1,\ldots,n_r}$ as a rational map onto its image,
then we have that  $\mathfrak S_{n_1,\ldots,n_r}$ is an isomorphism.
Indeed, for $j=1,\ldots, r$, the module of representatives of the $j$-th component $\mathfrak T_j$ of the inverse 
$\mathfrak T = \mathfrak S_{n_1,\ldots,n_r}^{-1}$ is 
generated by the $(n_1+1)\cdots(n_{j-1}+1) (n_{j+1}+1)\cdots (n_r+1)$ vectors
$(z_{(\iota_1,\ldots,\iota_r)}:\iota_j=0,\ldots,n_j)$,
as $\iota_1,\ldots,\iota_{j-1},\iota_{j+1},\ldots,\iota_{r}$ vary. Note, in particular, that
$\mathfrak T_j$ is not uniquely represented up to proportionality, 
provided that $n_1,\ldots,n_{j-1},n_{j+1},\ldots,n_r$ are not all zero.

\subsubsection{Multi-rational maps as ordinary rational maps}\label{multi vs ordinary}
Let $\Phi=(\Phi_1,\ldots,\Phi_s):X\dashrightarrow Y$ be a multi-rational map. Then, by composing $\Phi$
with the restriction to $Y$ of the 
Segre embedding $
 \mathfrak S_{m_1,\ldots,m_s}:\PP^{m_1}\times \cdots \times \PP^{m_s}\longrightarrow \PP^M
 $, where 
 $M=(m_1+1)\cdots(m_s+1)-1$,    
  we get an ordinary rational map 
  $\widetilde{\Phi}:X\dashrightarrow \mathfrak S_{m_1,\ldots,m_s}(Y)\subseteq\PP^M$.
  The rational map 
 $\widetilde{\Phi}$ is the unique rational map 
   that makes the following diagram commutative:
 \begin{equation*} 
\xymatrix{ X \ar@{-->}[rrr]^{\widetilde{\Phi}} \ar@{-->}[rrd]^{\Phi_1} \ar@{.>}[rrrd] \ar@{-->}[rrrrd]^{\Phi_s}& & & \mathfrak{S}_{m_1,\ldots,m_s}(\PP^{m_1}\times \cdots \times \PP^{m_s})   \ar@/^/[ld] \ar@{.>}[d] \ar@/_/[rd] &  \\ 
 & & \PP^{m_1} & \cdots & \PP^{m_s}
} 
\end{equation*}
Since $\mathfrak S_{m_1,\ldots,m_s}$
 is an isomorphism onto its image, we have that $\Phi$ is a morphism (resp.,  birational; resp., isomorphism)
 if and only if $\widetilde{\Phi}$ is a morphism (resp., birational; resp., isomorphism).
 Thus,
 from a theoretical point of view, 
 it would be enough to consider only 
 ``\emph{ordinary}'' rational maps.
In practice, however, this complicates things considerably since
the ambient space of the target of $\mathfrak S_{m_1,\ldots,m_s}$ is much larger with respect to the source, and moreover 
the homogeneous coordinate ring 
of the image of $\mathfrak S_{m_1,\ldots,m_s}$ is no longer more a
unique factorization domain (ruling out trivial cases).

\subsubsection{Graph of a (multi)-rational map}\label{subsect graph}
Let $\mathbf F^{(1)},\ldots,\mathbf F^{(s)}$ be, respectively, representatives 
of the components $\Phi_1,\ldots,\Phi_s$ of a multi-rational map $\Phi:X\dashrightarrow Y$.
Consider the $\ZZ^r\times \ZZ^s$-graded coordinate ring of 
\begin{equation}\label{PPP}
\PP^{n_1}\times \cdots \times \PP^{n_r}\times\PP^{m_1}\times \cdots \times \PP^{m_s},
\end{equation} 
given by
\[T =  K[\mathbf{x}_1;\ldots;\mathbf{x}_r;\mathbf{y}_1;\ldots;\mathbf{y}_s],\]
where $\mathbf{x}_{j} = (x_0^{(j)},\ldots,x_{n_j}^{(j)})$
and $\mathbf{y}_{i} = (y_0^{(i)},\ldots,y_{m_i}^{(i)})$, for $j=1,\ldots,r$ and $i=1,\ldots,s$.
Moreover, let $t_1,\ldots,t_s$ be new variables, and consider 
the extended polynomial ring 
\[\overline{T} =  K[t_1,\ldots,t_s;\mathbf{x}_1;\ldots;\mathbf{x}_r;\mathbf{y}_1;\ldots;\mathbf{y}_s] . \]
We define an ideal in $\overline{T}$ as the following sum of ideals
(by abuse of notation we denote by $\mathbf F^{(i)}$ also chosen lifts of $\mathbf F^{(i)}$ to $R$):
\begin{equation}\label{PPee}
 \mathcal{I}_{(\mathbf{F}^{(1)},\ldots,\mathbf{F}^{(s)})} := I + \sum_{i=1}^s \left(\mathbf{y}_i -t_i\,\mathbf{F}^{(i)}\right).
\end{equation}
The \emph{graph} $\Gamma(\Phi)$ of the multi-rational map $\Phi$ is the subvariety of \eqref{PPP}
defined by the contraction ideal 
\begin{equation}
 \mathcal{I}_{(\mathbf{F}^{(1)},\ldots,\mathbf{F}^{(s)})}\cap T , 
\end{equation}
which no longer depends on the choice of the representatives $\mathbf{F}^{(i)}$.
Equivalently, we can consider the homogeneous ideal in $T$ given by 
\begin{equation}\label{graphSat1}
 \mathcal{J}_{(\mathbf{F}^{(1)},\ldots,\mathbf{F}^{(s)})} := I + \left(\mbox{$2\times 2$ minors of } \begin{pmatrix} y_0^{(i)} & \cdots & y_{m_i}^{(i)} \\ F_0^{(i)} & \cdots & F_{m_i}^{(i)} \end{pmatrix},\ i=1,\ldots,s \right) ,
 \end{equation}
and then we can calculate the ideal of  $\Gamma(\Phi)$ by the saturation:
 \begin{equation}\label{graphSat2}
 \left(\cdots\left( \mathcal{J}_{(\mathbf{F}^{(1)},\ldots,\mathbf{F}^{(s)})} : (\mathbf{F}^{(1)})^{\infty}\right):\cdots\right): (\mathbf{F}^{(s)})^{\infty} .
 \end{equation}
 We point out that the homogeneous coordinate ring of $\Gamma(\Phi)$
 is also known as ``\emph{Rees algebra}'', see \cite{ReesIdeal} and references therein.
    
 We have two projections (which are morphisms) that fit in a commutative diagram
\begin{equation*} 
\xymatrix{ & \Gamma(\Phi) \ar[dl]_{\pi_1} \ar[dr]^{\pi_2}\\ X \ar@{-->}[rr]^{\Phi}& & Y } 
\end{equation*}
The first projection $\pi_1:\Gamma(\Phi)\to X$ 
  is also known as the \emph{blowing up of $X$ along $B$},
  where $B=X\setminus\Dom(\Phi)$ is the base locus of $\Phi$.
  It is a birational morphism,
  and it is an isomorphism if and only if $\Phi$ is a morphism.
  See \emph{e.g.} \cite[Chapter~II, Section~7]{hartshorne-ag} for more details.
 The second projection  $\pi_2:\Gamma(\Phi)\to Y$ is birational if and only if $\Phi$ is birational,
 and in that case 
 the graph of $\Phi^{-1}$ is the same as
  that of $\Phi$, by exchanging the two projections.
 Moreover,
$\pi_2$ and $\Phi$ have always the same image in $Y$; 
   in particular, we can calculate the homogeneous ideal  of the image of $\Phi$ as
   the
   contraction of the ideal of $\Gamma(\Phi)$ to $S=K[\mathbf{y}_1;\ldots;\mathbf{y}_s]$.

 \subsubsection{Computing the inverse map of a birational map}
 Keep the notation as above, and assume moreover that 
 $\Phi:X\dashrightarrow Y$ is birational. 
 We want to find the components $\Psi_j:Y\dashrightarrow \mathbb{P}^{n_j}$, for $j=1,\ldots,r$,
 of the inverse multi-rational map $\Psi:Y\dashrightarrow X$ of $\Phi$.
  
 Fix a minimal set of multi-forms 
 generating  the homogeneous ideal of the graph $\Gamma(\Phi)$ 
 in 
 the $\ZZ^r\times \ZZ^s$-graded coordinate ring of \eqref{PPP}.
 For each $j=1,\ldots,r$, we select in this set those of 
 multi-degree $(0,\ldots,0,1,0,\ldots,0;d_1,\ldots,d_s)$,
 where $1$ occurs at position $j$, and $d_1,\ldots,d_s$ are not subject to conditions.
 Let us denote 
 these multi-forms by
  $H_1(\mathbf x_j,\mathbf y_1,\ldots,\mathbf y_s),\ldots,H_q(\mathbf x_j,\mathbf y_1,\ldots,\mathbf y_s)$.
 Thus, for $k=1,\ldots,q$, we can write  
 \[H_k(\mathbf x_j,\mathbf y_1,\ldots,\mathbf y_s)=x_0^{(j)} {G}^{(j,k)}_0(\mathbf y_1,\ldots,\mathbf y_s)+\cdots+x_{n_j}^{(j)} G_{n_j}^{(j,k)}(\mathbf y_1,\ldots,\mathbf y_s),\]
 for suitable uniquely determined forms $G_{\iota_j}^{(j,k)}\in S=K[\mathbf y_1,\ldots,\mathbf y_s]$.
 We regard 
 the $q\times(n_j+1)$-matrix 
 \[{\mathfrak J}^{(j)} = \left(G_{\iota_j}^{(j,k)}\right)^{\iota_j=0,\ldots,n_j}_{k=1,\ldots,q}\]
 as a 
 matrix over the 
 homogeneous coordinate ring $S/J$ of $Y$.
 We have the following: 
 \begin{proposition}
  The $(S/J)$-module of representatives of 
  $\Psi_j
  $
  is given by $\ker({\mathfrak J}^{(j)})$. 
  More explicitly we have that the
rank of 
${\mathfrak J}^{(j)}$ is $n_j$, and
 $\Psi_j$ is represented by 
 the vector of signed $n_j\times n_j$-minors
 of any full rank $n_j\times(n_j+1)$-submatrix of ${\mathfrak J}^{(j)}$.
 \end{proposition}
 A proof of the previous result can be found 
in \cite[Theorem~2.4]{Simis2004162},
in the particular case when $r=s=1$ (see also \cite{DORIA2012390} and \cite[Theorem~4.4]{BusY2020} for the case when $s=1$ and the source is a product of projective varieties).
The proof in the general case is not so different;
its main ingredients are: 
the description of the equations of the graph $\Gamma(\Phi)$  given by 
\eqref{graphSat1} and \eqref{graphSat2}, and the fact 
that  $\Gamma(\Phi)$
can be identified with $\Gamma(\Psi)$.
We leave the details to the reader.

 \subsubsection{Direct and inverse images via multi-rational maps}
  If $Z\subseteq X$ is an irreducible subvariety such that $Z\cap\Dom(\Phi)\neq \emptyset$,
 then we can consider the {restriction of $\Phi$ to $Z$}, $\Phi|_Z:Z\dashrightarrow Y$,
 defined as usual by the composition 
of the inclusion $Z\hookrightarrow X$ with $\Phi$.
 Note that 
the graph (and hence the image) of $\Phi|_Z$, can be calculated as above,
just by replacing in \eqref{PPee} the ideal  $I$  with the multi-saturated homogeneous ideal of  $Z$, 
and by choosing the representatives $\mathbf{F}^{(i)}$ such that
$Z \nsubseteq  V(\mathbf{F}^{(i)})$.
This gives us a way to calculate the direct image $\overline{\Phi(Z)}=\overline{\Phi|_Z(Z)}$.

If $W\subseteq Y$ is a subvariety,
we can calculate the inverse image $\overline{\Phi^{-1}(W)}\subseteq  X$ 
as $\overline{\Phi^{-1}(W)}=\overline{\widetilde{\Phi}^{-1}(\mathfrak{S}_{m_1,\ldots,m_s}(W))}$,
using  Proposition~\ref{direct and inverse image}.
Alternatively (and more efficiently), let 
$I_W\subseteq S$ be the defining ideal of $W$,
and let 
$\varphi_{(\mathbf F^{(1)},\ldots,\mathbf F^{(s)})}:S\to R/I$ be the map 
 defined by $y_{\iota_i}^{(i)} \mapsto F_{\iota_i}^{(i)}\in R/I$,
 for $i=1,\ldots,s$ and $\iota_i=0,\ldots,m_i$.
 Then the saturation of
 the extended ideal $\left(\varphi_{(\mathbf F^{(1)},\ldots,\mathbf F^{(s)})}(I_W)\right)\subseteq R/I$
 with respect to all the ideals 
 $(\mathbf F^{(i)})$,
 for $i = 1,\ldots,s$, gives us the ideal of 
 the closure of
 ${\overline{\Phi^{-1}(W)}\setminus V(\mathbf F^{(1)},\ldots,\mathbf F^{(s)})}$. 

  \subsubsection{Multi-degree of a multi-rational map}
  Let $\Phi:X\dashrightarrow Y$ be a rational map.
 The \emph{projective degrees} $d_0(\Phi),d_1(\Phi),\ldots,d_{\dim X}(\Phi)$ of $\Phi$ are 
  defined as the components of 
  the multi-degree of the graph, 
  embedded as a subvariety 
  of 
  \[\mathfrak{S}_{n_1,\ldots,n_r}(\mathbb{P}^{n_1}\times\cdots\times \mathbb{P}^{n_r})\times
  \mathfrak{S}_{m_1,\ldots,m_s}(\mathbb{P}^{m_1}\times\cdots\times \mathbb{P}^{m_s})
  \subset \PP^{N}\times \PP^{M} ,\]
  where $N=\Pi_{j=1}^r (n_j+1) -1$ and $M=\Pi_{i=1}^s (m_i+1) -1$.
It follows that the composition 
$\widetilde{\Phi}:X\dashrightarrow\PP^{M}$ 
of $\Phi$
with the restriction to $Y$ of the Segre embedding $\mathfrak{S}_{m_1,\ldots,m_s}$
has the same projective degrees as $\Phi$.
If $L$ denotes the intersection of $Y$ with $\dim X - i$
general hypersurfaces of multi-degree $(1,\ldots,1)$, then we have
\[
 d_i(\Phi)= \deg \left(\mathfrak{S}_{n_1,\ldots,n_r}(\overline{\Phi^{-1}(L)})\right),
\]
if $\dim \left(\overline{\Phi^{-1}(L)}\right) = i $ and 
 $d_i(\Phi)=0$ otherwise.
See also \cite[Example~19.4, p.~240]{harris-firstcourse}.
 This gives us a probabilistic algorithm to compute the projective degrees, as already remarked in \cite{packageCremona}.
A non-probabilistic algorithm  can be obtained by calculating the multi-degree 
of the graph of $\Phi$ as a subvariety of $\mathbf P^{n_1,\ldots,n_r}\times\mathbf P^{m_1,\ldots,m_s}$ 
and then applying the following remark.
\begin{remark}
  Let $P(a_1,\ldots,a_r,b_1,\ldots,b_s)\in\mathbb{Z}[a_1,\ldots,a_r,b_1,\ldots,b_s]$ 
  be the multi-degree of a $k$-dimensional subvariety of 
  $\mathbf P^{n_1,\ldots,n_r}\times \mathbf P^{m_1,\ldots,m_s}$.
  Then the multi-degree of the same variety 
  embedded as a subvariety of 
  $\mathfrak{S}_{n_1,\ldots,n_r}(\mathbf P^{n_1,\ldots,n_r})\times\mathfrak{S}_{m_1,\ldots,m_s}(\mathbf P^{m_1,\ldots,m_s})\subset \PP^{N}\times \PP^M$, is given by
  \[\sum_{i=\max(0,k-M)}^{\min(k,N)} d_i \, a^{N-i}  b^{M-k+i}\in \mathbb{Z}[a,b],\]
  where
  $d_i$ denotes the coefficient of the monomial $a_1^{n_1}\cdots a_r^{n_r} b_1^{m_1}\cdots b_s^{m_s}$
  in the polynomial $(a_1+\cdots+a_r)^i (b_1+\cdots+b_s)^{k -i} P(a_1,\ldots,a_r,b_1,\ldots,b_s)$.
 In particular, when $m_1=\cdots=m_s=0$ we get 
  the degree 
   of the variety embedded in $\PP^N$ 
  from its multi-degree as a subvariety of $\mathbf P^{n_1,\ldots,n_r}$.
 \end{remark}
 The last projective degree  $d_{\dim X}(\Phi)$
  is the degree of $\mathfrak S_{n_1,\ldots,n_r}(X)\subseteq\PP^N$.
  The first projective degree $d_0(\Phi)$ is the product of the degree of 
  $\mathfrak S_{m_1,\ldots,m_s}(\overline{\Phi(X)})\subseteq\PP^M$
  with the \emph{degree} of $\Phi$.
  We have that $\Phi$ is birational onto its image 
  if and only if its 
  degree is $1$, that is, if and only if
  $d_0(\Phi) = \deg\left(\mathfrak S_{m_1,\ldots,m_s}(\overline{\Phi(X)})\right)$.
  Thus we can determine whether  $\Phi$ is birational without computing its inverse.

  \section{Implementation in {\it Macaulay2}}\label{section 2}
  The  \emph{Macaulay2} package 
  \href{https://faculty.math.illinois.edu/Macaulay2/doc/Macaulay2/share/doc/Macaulay2/MultiprojectiveVarieties/html/index.html}{\emph{MultiprojectiveVarieties}}
  provides support for multi-projective varieties and multi-rational maps.
  It implements, among other things, the methods described in the previous section.
 As we previously said, a multi-rational map can be represented by a list of rational maps 
  having as target a projective space. 
  Partial support for this particular kind of rational
  maps
  is provided by the package 
  \href{https://faculty.math.illinois.edu/Macaulay2/doc/Macaulay2/share/doc/Macaulay2/Cremona/html/index.html}{\emph{Cremona}} \cite{CremonaPackageSource},
  on which the first one depends.
  
  Here we give just one simple example to illustrate
  how one can work with these packages. We refer to the online documentation of \emph{Macaulay2} 
  for more examples and technical details.
\medskip 

  It is classically well known that 
  a smooth cubic hypersurface $X\subset\PP^5$
  containing two disjoint planes is birational to $\PP^2\times\PP^2$,
  and that the inverse map $\PP^2\times\PP^2\dashrightarrow X$ 
  is not defined along a K3 surface of degree $14$. We now analyze this example using \emph{Macaulay2}.

  In the following lines of code, we first define the two projections $f:\PP^5\dashrightarrow\PP^2$
  and $g:\PP^5\dashrightarrow\PP^2$ 
  from two disjoint planes in $\PP^5$, then we define 
  the multi-rational map $(f,g):\PP^5\dashrightarrow\PP^2\times\PP^2$ 
  and restrict it to a smooth cubic hypersurface $X$
  containing the two planes. So we get a multi-rational map $\Phi:X\dashrightarrow\PP^2\times\PP^2$.
   \begin{tcolorbox}[breakable=true,boxrule=0pt,opacityback=0.0,enhanced jigsaw]
{\footnotesize
\begin{Verbatim}[commandchars=&!$]
&colore!darkorange$!M2 --no-preload$
&colore!output$!Macaulay2, version 1.18$
&colore!darkorange$!i1 :$ &colore!airforceblue$!needsPackage$ "&colore!bleudefrance$!MultiprojectiveVarieties$"; &colore!commentoutput$!-- version 2.2$
&colore!darkorange$!i2 :$ K = &colore!airforceblue$!QQ$, K[t,u,v,x,y,z];
&colore!darkorange$!i3 :$ f = &colore!bleudefrance$!rationalMap$ {t,u,v};
&colore!circOut$!o3 :$ &colore!output$!RationalMap (linear rational map from PP^5 to PP^2)$
&colore!darkorange$!i4 :$ g = &colore!bleudefrance$!rationalMap$ {x,y,z};
&colore!circOut$!o4 :$ &colore!output$!RationalMap (linear rational map from PP^5 to PP^2)$
&colore!darkorange$!i5 :$ Phi = &colore!bleudefrance$!rationalMap$ {f,g};
&colore!circOut$!o5 :$ &colore!output$!MultirationalMap (rational map from PP^5 to PP^2 x PP^2)$
&colore!darkorange$!i6 :$ X = &colore!bleudefrance$!projectiveVariety$ &colore!airforceblue$!ideal$(t*u*x-u^2*x+u*v*x-v^2*x+t*x^2-u*x^2+t^2*y-t*u*y-
               t*v*y-t*x*y-v*x*y-t*y^2+t*u*z+v^2*z-t*x*z-u*y*z-v*y*z-t*z^2+u*z^2);
&colore!circOut$!o6 :$ &colore!output$!ProjectiveVariety, hypersurface in PP^5$
&colore!darkorange$!i7 :$ Phi = Phi|X;
&colore!circOut$!o7 :$ &colore!output$!MultirationalMap (rational map from X to PP^2 x PP^2)$
\end{Verbatim}
} 
\end{tcolorbox}\noindent 
Next, we verify that $\Phi$ is dominant and birational, compute the inverse map $\Phi^{-1}$,
and ``describe'' the base locus of $\Phi^{-1}$.
 \begin{tcolorbox}[breakable=true,boxrule=0pt,opacityback=0.0,enhanced jigsaw]
{\footnotesize
\begin{Verbatim}[commandchars=&!$]
&colore!darkorange$!i8 :$ &colore!bleudefrance$!image$ Phi == &colore!bleudefrance$!target$ Phi
&colore!circOut$!o8 =$ &colore!output$!true$
&colore!darkorange$!i9 :$ &colore!bleudefrance$!degree$ Phi
&colore!circOut$!o9 =$ &colore!output$!1$
&colore!darkorange$!i10 :$ &colore!bleudefrance$!inverse$ Phi;
&colore!circOut$!o10 :$ &colore!output$!MultirationalMap (birational map from PP^2 x PP^2 to X)$
&colore!darkorange$!i11 :$ &colore!bleudefrance$!describe baseLocus inverse$ Phi;
&colore!circOut$!o11 =$ &colore!output$!ambient:.............. PP^2 x PP^2$
&colore!output$!      dim:.................. 2$
&colore!output$!      codim:................ 2$
&colore!output$!      degree:............... 14$
&colore!output$!      multidegree:.......... 2 T_0^2 + 5 T_0 T_1 + 2 T_1^2$
&colore!output$!      generators:........... (2,1)^1 (1,2)^1$
&colore!output$!      purity:............... true$
&colore!output$!      dim sing. l.:......... -1$
\end{Verbatim}
} 
\end{tcolorbox}\noindent 
Now we take the graph of $\Phi$ with the two projections 
$p_1:\Gamma(\Phi)\to X$ and $p_2:\Gamma(\Phi)\to\PP^2\times\PP^2$. We calculate
the projective degrees of $p_1$ and $p_2$, the inverse of $p_2$,
and verify that $p_1\circ p_2^{-1} = \Phi^{-1}$ and that $p_2$ is a morphism but not an isomorphism.
 \begin{tcolorbox}[breakable=true,boxrule=0pt,opacityback=0.0,enhanced jigsaw]
{\footnotesize
\begin{Verbatim}[commandchars=&!$]
&colore!darkorange$!i12 :$ (p1,p2) = &colore!bleudefrance$!graph$ Phi;
&colore!darkorange$!i13 :$ (&colore!bleudefrance$!multidegree$ p1, &colore!bleudefrance$!multidegree$ p2)
&colore!circOut$!o13 =$ &colore!output$!({141, 63, 25, 9, 3}, {141, 78, 40, 18, 6})$
&colore!darkorange$!i14 :$ &colore!bleudefrance$!inverse$ p2;
&colore!circOut$!o14 :$ &colore!output$!MultirationalMap (birational map from PP^2 x PP^2 to 4-dimensional$
&colore!output$!                        subvariety of PP^5 x PP^2 x PP^2)$
&colore!darkorange$!i15 :$ (&colore!bleudefrance$!inverse$ p2) * p1 == &colore!bleudefrance$!inverse$ Phi, &colore!bleudefrance$!isMorphism$ p2, &colore!bleudefrance$!isIsomorphism$ p2
&colore!circOut$!o15 =$ &colore!output$!(true, true, false)$
\end{Verbatim}
} 
\end{tcolorbox}\noindent 
We now calculate the \emph{exceptional locus} of the first projection $p_1$;
this is the inverse image of the base locus of $p_1^{-1}$. 
 \begin{tcolorbox}[breakable=true,boxrule=0pt,opacityback=0.0,enhanced jigsaw]
{\footnotesize
\begin{Verbatim}[commandchars=&!$]
&colore!darkorange$!i16 :$ &colore!bleudefrance$!baseLocus$ Phi == &colore!bleudefrance$!baseLocus$ &colore!bleudefrance$!inverse$ p1
&colore!circOut$!o16 =$ &colore!output$!true$
&colore!darkorange$!i17 :$ E = p1^* (&colore!bleudefrance$!baseLocus$ Phi);
&colore!circOut$!o17 :$ &colore!output$!ProjectiveVariety, threefold in PP^5 x PP^2 x PP^2$
&colore!darkorange$!i18 :$ &colore!bleudefrance$!dim$ E, &colore!bleudefrance$!degree$ E
&colore!circOut$!o18 =$ &colore!output$!(3, 48)$
\end{Verbatim}
} 
\end{tcolorbox}\noindent 
Finally, we take the first projection $h:\Gamma(p_2)\to \Gamma(\Phi)$ from the graph of $p_2$.
This multi-rational map, regarded as a rational 
map between embedded projective varieties,
has as source a fourfold of degree $771$ in $\PP^{485}$
and as target a fourfold of degree $141$ in $\PP^{53}$.
 \begin{tcolorbox}[breakable=true,boxrule=0pt,opacityback=0.0,enhanced jigsaw]
{\footnotesize
\begin{Verbatim}[commandchars=&!$]
&colore!darkorange$!i19 :$ h = &colore!airforceblue$!first$ &colore!bleudefrance$!graph$ p2;
&colore!circOut$!o19 :$ &colore!output$!MultirationalMap (birational map from 4-dimensional subvariety of$
&colore!output$!                        PP^5 x PP^2 x PP^2 x PP^2 x PP^2 to 4-dimensional$
&colore!output$!                        subvariety of PP^5 x PP^2 x PP^2)$
&colore!darkorange$!i20 :$ &colore!bleudefrance$!degree source$ h, &colore!bleudefrance$!degree target$ h
&colore!circOut$!o20 =$ &colore!output$!(771, 141)$
\end{Verbatim}
} 
\end{tcolorbox}\noindent 
By construction, we know (and \emph{Macaulay2} knows) that the map $h$ is birational. We can also verify this experimentally,
by reducing to prime characteristic and calculating the fiber of $h$ at a random point $p$ on its source.
\begin{tcolorbox}[breakable=true,boxrule=0pt,opacityback=0.0,enhanced jigsaw]
{\footnotesize
\begin{Verbatim}[commandchars=&!$]
&colore!darkorange$!i21 :$ h = h ** (&colore!airforceblue$!ZZ$/1000003),;
&colore!darkorange$!i22 :$ p = &colore!bleudefrance$!point source$ h;
&colore!circOut$!o22 =$ &colore!output$!ProjectiveVariety, a point in PP^5 x PP^2 x PP^2 x PP^2 x PP^2$
&colore!darkorange$!i23 :$ p == h^* h p
&colore!circOut$!o23 =$ &colore!output$!true$
\end{Verbatim}
} 
\end{tcolorbox}\noindent 
On a standard laptop,
the time to execute the
23 lines of code above 
is less than 5 seconds.
 

\providecommand{\bysame}{\leavevmode\hbox to3em{\hrulefill}\thinspace}
\providecommand{\MR}{\relax\ifhmode\unskip\space\fi MR }
\providecommand{\MRhref}[2]{%
  \href{http://www.ams.org/mathscinet-getitem?mr=#1}{#2}
}
\providecommand{\href}[2]{#2}

\end{document}